\documentclass[conference]{IEEEtran}
\usepackage{amsmath,amsfonts,amssymb}
\usepackage{graphicx}
\usepackage{cite}
\usepackage{tikz}
\usepackage{algorithm}
\usepackage{algorithmic}
\title{Optimization of Bottlenecks in Quantum Graphs Guided by Fiedler Vector-Based Spectral Derivatives}

\author{
\IEEEauthorblockN{John T.M CAMPBELL, JOHN DOOLEY}
\IEEEauthorblockA{Maynooth Dept of Electronic Engineering \\
Email: \{John.Campbell.2023@mumail.com\}}
}

\begin{document}

\maketitle

\begin{abstract}
This paper discusses the relationships between the Fiedler vector, the Cheeger constant, and threshold behaviors in networks of quantum resource nodes represented as Quantum Directed Acyclic Graphs (QDAGs). We explore how these mathematical constructs can be applied to understand the dynamics of quantum information flow in QDAGs, especially in the context of routing problems with bottlenecks in graph signal processing, and how new eigenvalue-based rewiring techniques can optimize entanglement distribution between nodes in a QDAG.
\end{abstract}

\IEEEpeerreviewmaketitle

\section{Introduction}
Quantum Directed Acyclic Graphs (QDAGs) are a versatile framework in entanglement networking that encode relationships between entanglement resources, nonlinear oscillators elements, i.e., qubit sources, memories, repeaters, clients etc., in a non-cyclic manner \cite{2}. Understanding the underlying mathematical structures governing QDAGs allows for optimization using algorithms developed from classical graph theory and network science \cite{vanderberg2012} to help understand and optimize quantum information flow. Central to this discussion are the Fiedler vector and the Cheeger constant [5], which provide insight into the connectivity and bottlenecks in a network, two details which may not always be mutually exclusive as this research hopes to elucidate. A research focus is to examine implications for quantum states on the graph i and when the Fiedler eigenbasis is changed. This involves examining dynamical aspects of quantum systems represented by graphs, properties such as stability, connectivity as they are rewired for optimization .To accomplish this we start at the foundations of the Fiedler vector, in theorem form, to form conceptual bridges with the field of dynamical quantum graphs. Also demonstrated is the use of the Fiedler vector to perform a novel type of graph rewiring for optimization of Quantum DAG's. This method is dubbed QFERN, standing for Quantum Fiedler Eigenbasis for Rewiring Networks. It represents a method, implimented as a workflow and then as an algorithm, that utilizes specifically the Fiedler vector and eigenbasis, to rewire or reconfigure quantum networks to optimize the bottlenecks in such networks, as represented by spectral gaps, and effective resistances in such graphs of the entanglement between nodes.  

\section{Mathematical Foundations}

\subsection{Fiedler Vector}
The Fiedler vector \(\mathbf{v}\) is derived from the second smallest eigenvector of the Laplacian matrix \(\mathbf{L}\) of a graph. The Laplacian is given by:

\[
\mathbf{L} = \mathbf{D} - \mathbf{A}
\]

where \(\mathbf{D}\) is the degree matrix, and \(\mathbf{A}\) is the adjacency matrix. The Fiedler vector captures essential connectivity properties of the graph, indicating how tightly coupled the vertices are \cite{2}.

\subsection{Cheeger Constant}
The Cheeger constant \(h\) quantifies the connectivity of a graph and is defined as:

\[
h = \min_{S \subset V} \frac{|\partial S|}{|S|}
\]

where \(|\partial S|\) counts the edges crossing the cut defined by \(S\) and \(V\) is the set of vertices. The Cheeger constant is crucial for understanding the bottleneck behaviors in QDAGs \cite{cheeger}.

\section{Relation between Fiedler Vector and Cheeger Constant}
The Fiedler vector and Cheeger constant are related through spectral graph theory. Specifically, the Cheeger inequality provides a connection:

\[
\frac{h^2}{2} \leq \lambda_2 \leq 2h
\]

where \(\lambda_2\) is the second smallest eigenvalue of the Laplacian matrix. A larger Cheeger constant reflects a better connectivity, aligning with a lower value of the Fiedler vector's components \cite{vanderberg2012}.

\subsection{Fiedler Vector and Dirichlet Energies}

The Fiedler vector, denoted as \(\mathbf{f}_2 \in \mathbb{R}^n\), corresponds to the second smallest eigenvector of the Laplacian matrix of a graph, often termed the Fiedler matrix. This vector provides essential insights into the connectivity properties of the graph, particularly in terms of cluster membership \cite{vanderberg2012}. Specifically, the elements of \(\mathbf{f}_2\) can be interpreted as indicators of node membership to each of the two clusters that arise from a spectral clustering process.

The associated eigenvalue \(\lambda_2\) of \(\mathbf{f}_2\) encapsulates the Dirichlet energy of the system, providing a measure of the smoothness of the function represented by the Fiedler vector across the graph. This eigenvalue plays a critical role in quantifying the overall "bottleneckedness" of the graph, as described by the Cheeger constant, which, in turn, reflects the minimum ratio of the edge cut to the volume of the smaller subset of vertices within the graph.

We aim to optimize the bottleneck width through a novel form of diffusion-based graph rewiring characterized by a search matrix \(\mathbf{A}_{\text{soft}}\), which closely resembles the original adjacency matrix \(\mathbf{A}\) but minimizes the bottleneck size. This optimization is guided by spectral derivatives, allowing us to express the adjustments in the graph structure as follows:

\[
\mathcal{L}_{\text{Fiedler}} = \mathbf{A}_{\text{soft}} - \mathbf{A} + \alpha \frac{\lambda_2}{2} \nabla \mathbf{A}_{\text{soft}}\lambda_2 \equiv \text{Tr} \nabla \mathcal{L}_{\text{mi}} \lambda_2
\]

To elucidate the relationship between the Fiedler vector and the structure of the graph, we employ the relationship:

\[
T \nabla \mathbf{A}_{\text{soft}} \mathcal{L}_{\text{mi}} = \text{diag}(\mathbf{f}_2 \mathbf{f}_2^T) \mathbf{1}_{1}^T - \mathbf{f}_2 \mathbf{f}_2^T
\]

The expression for \(\lambda_2\) can thus be seen as representing the spectral gap or bottleneck size, indicating how effectively the clusters are separated in the process of optimization. 

As we manipulate the underlying graph structure, the evolution of the Fiedler vector and its corresponding eigenvalue \(\lambda_2\) plays a crucial role in guiding the optimization of the graph to achieve a balance between connectivity and separability of the clusters present.

A proposed algorithm to accomplish this requires the following, general recipe for graphs: 

\subsection{Graph Processing Steps}
\begin{itemize}
    \item \textbf{Random Directed Acyclic Graph Generation:} We use NetworkX to create a random directed acyclic graph (DAG).
    \item \textbf{Fiedler Vector Calculation:} We compute the Fiedler vector from the Laplacian matrix of the original adjacency matrix. The Fiedler vector is the second smallest eigenvector of the Laplacian, which is critical for our adjustments.
    \item \textbf{Gradient Calculation:} The gradient is derived based on the outer product of the Fiedler vector (-$|fiedler\_vector[i] * fiedler\_vector[j]|$). This allows us to adjust the adjacency matrix in a manner reflecting the graph's structural properties.
    \item \textbf{Adjustment of Asoft:} The Asoft matrix is adjusted using the gradient while maintaining non-negativity and symmetry. The adjustment is controlled by the alpha parameter.
    \item \textbf{Laplacian Calculation and Plotting:} The code calculates the Laplacian for both the original adjacency matrix and the Asoft, then plots them for visual comparison  as in Fig 1.
\end{itemize}

\begin{figure}[ht]
    \centering
    \includegraphics[width=0.50\textwidth]{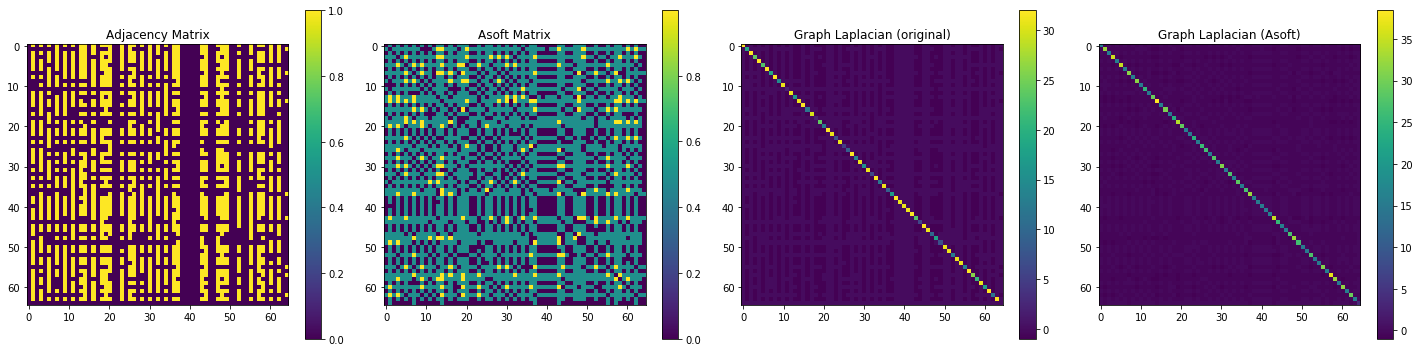}
    \caption{Graph Adjacency and Asoft Matrices and their respective Graph Laplacians}
    \label{fig:unit_circle_pulses}
\end{figure}

\section{Quantum Directed Acyclic Graphs (QDAGs)}
QDAGs are employed to model quantum circuits. The vertices represent quantum gates, and the edges signify the quantum states’ directed flow. The Fiedler vector evaluates the potential bottlenecks in these flows \cite{2}.

\subsection{Thresholds in QDAGs}
In the context of QDAGs, we define critical bottlenecks as thresholds for qubit entanglement and coherence loss. These thresholds can be derived from the Cheeger constant, treating it as a bottleneck index, thus linking classical graph properties with quantum behaviors, such as bottlenecks leading to effective resistances to the transfer of entanglement in a network.

\subsection{Effective Resistances in QDAGs}

In the context of Quantum Directional Acyclic Graphs (QDAGs), the effective resistance between two nodes \( u \) and \( v \) can be quantitatively analyzed using the properties of the graph's Laplacian matrix. The effective resistance \( R_{u,v} \) is mathematically expressed as:

\begin{equation}
    R_{u,v} = (e_u - e_v)^T L^+ (e_u - e_v),
\end{equation}

where \( L^+ \) denotes the pseudoinverse of the Laplacian matrix \( L \). This formulation emphasizes the importance of the Laplacian in characterizing the resistance in a network by capturing the influence of the graph's topology on the flow between nodes.

To deepen our understanding of effective resistance in the framework of QDAGs, we can invoke the relationship between effective resistance and the spectral characteristics of the graph. In particular, we can connect the effective resistance to the eigenvalues associated with the Fiedler vector, which provides insight into the graph's connectivity. This relationship can be formulated as:

\begin{equation}
    R_{u,v} = \sigma_i > 2 \sum_{i=1}^{n} \frac{1}{\lambda_i} f_i f_i^T,
\end{equation}

In this equation, \( \sigma_i \) represents the weights assigned to the eigenvalues \( \lambda_i \), and \( f_i \) corresponds to the components of the Fiedler vector tied to the nodes of the graph. The Fiedler vector, being the eigenvector associated with the second smallest eigenvalue of the Laplacian, plays a pivotal role in describing the connectivity of the network \cite{2}.

This connection between the structure of QDAGs, their spectral properties, and the effective resistance highlights how the Fiedler vector encapsulates the underlying entanglement flow characteristics within the network.

Optimizing the network Laplacian based on the Fiedler vector can lead to more robust designs of quantum networks that minimize error rates and maximize entanglement distribution in bottlenecked networks, i.e., with only a single entanglement resource and many clients, repeaters, memory nodes, etc.

\section{QFERN Method for Graph Modification wrt bottlenecks and effective resistances}

\subsection{Overview}
Given a directed acyclic graph (DAG), we aim to modify its adjacency matrix based on the spectral properties of its Laplacian matrix to produce a new matrix, termed \( A_{\text{soft}} \). The algorithm sequentially computes gradients based on the Fiedler vector and applies them to adjust \( A_{\text{soft}} \).

\begin{algorithm}
\caption{Modification of Adjacency Matrix}
Generate a random DAG with $$ n $$ nodes.
 Compute the adjacency matrix $$ A $$ of the graph.
Compute the Laplacian matrix $$ L = D - A $$.
Calculate the eigenvalues and eigenvectors of $$ L $$.
Identify the Fiedler vector from the eigenvectors.
Initialize $$ A_{\text{soft}} = A $$.
For{each $$ (i, j) $$ in the nodes}
    If{$i \neq j$}
        Compute gradient based on Fiedler vector: $$ \text{gradient}[i,j] = -f[i] \cdot f[j] $$
    EndIf
EndFor
Update $$ A_{\text{soft}} \leftarrow A_{\text{soft}} + \alpha \cdot \text{gradient} $$
Ensure $$ A_{\text{soft}} $$ is non-negative and symmetric.
\end{algorithm}

\subsection{Optimization Scheme}

The optimization scheme aims to minimize the effective resistance across the network while maximizing the Cheeger constant, reflecting the network's bottlenecks. In this subsection, we will outline the procedure for this optimization process, including the necessary calculations and evaluations.

\subsubsection{Objective Function}

The optimization problem can be formulated as:

\begin{equation}
\min_{G} R_{\text{eff}}(G) \quad \text{s.t.} \quad \max_{G} h(G)
\end{equation}

where \( R_{\text{eff}} \) is the effective resistance of the graph \( G \), and \( h(G) \) is the Cheeger constant of the graph.

\subsubsection{Cheeger Constant Calculation}

The Cheeger constant \( h(G) \) is defined as:

\begin{equation}
h(G) = \min_{S \subseteq V} \frac{|\partial S|}{\min\{|S|, |V \setminus S|\}}
\end{equation}

where \( V \) is the set of vertices in the graph, \( |S| \) is the size of the set \( S \), and \( |\partial S| \) is the number of edges crossing the boundary of \( S \) (i.e., edges that have one endpoint in \( S \) and the other in \( V \setminus S \)).

The initial Cheeger constant is calculated to assess the connectivity of the network before any modifications are made.

\subsubsection{Optimization Loop}

The optimization scheme operates in an iterative loop as follows:

\begin{enumerate}
    \item A random edge \( e \) is removed from the graph \( G \).
    The edges are added to the graph using the function:
 the set of edges connecting nodes in groups A and B.
    \item A new edge \( e' \) is added from a list of possible edges, ensuring adherence to the graph's structural rules. Let \( E' \) be the set of possible edges. \begin{equation}
    E = \{ (i, j) | i \in A, j \in B \}
\end{equation}
    \item The new Cheeger constant \( h'(G) \) is calculated as per Equation (2) and printed.
\end{enumerate}

This loop continues until a stopping criterion is met, such as a maximum number of iterations or no improvement observed in the Cheeger constant.

\subsubsection{Objective Evaluation}

During the optimization process, we dynamically adjust the edges and evaluate configurations that yield better Cheeger constants while simultaneously analyzing the effective resistances:

\begin{equation}
R_{\text{eff}}(G) = \sum_{u, v \in V} R_{uv}
\end{equation}

where \( R_{uv} \) is the effective resistance between nodes \( u \) and \( v \).

The effective resistance $$R(u, v)$$ between two nodes $$u$$ and $$v$$ is calculated using the eigenvalues $$\lambda_i$$ and the Fiedler vector $$f$$ obtained from the Laplacian's eigenvalue decomposition:
\begin{equation}
    R(u, v) = \sum_{i=1}^{n} \frac{f_u f_v}{\lambda_i}
\end{equation}
where $$f_u$$ and $$f_v$$ are components of the Fiedler vector corresponding to nodes $$u$$ and $$v$$.

The optimization scheme examines how modifications in the edge structure can lead to favorable changes in both the Cheeger constant and effective resistance.

\subsubsection{Final Visualization}

Upon completing the QFERN optimization iterations, the effective resistances are recalculated to analyze the impact of edge modifications on network properties. Visualization techniques, such as graph plotting and heat maps, are employed to represent the changes in effective resistances and Cheeger constants across the iterations.

The resulting visualizations can help in understanding the effects of edge modifications and how they influence the overall connectivity and efficiency of the network.

We can compare the graph structure before and after the Adjacency to Asoft graph rewiring modifications to check for bottlenecks as follows:
\begin{figure}[ht]
    \centering
    \includegraphics[width=0.45\textwidth]{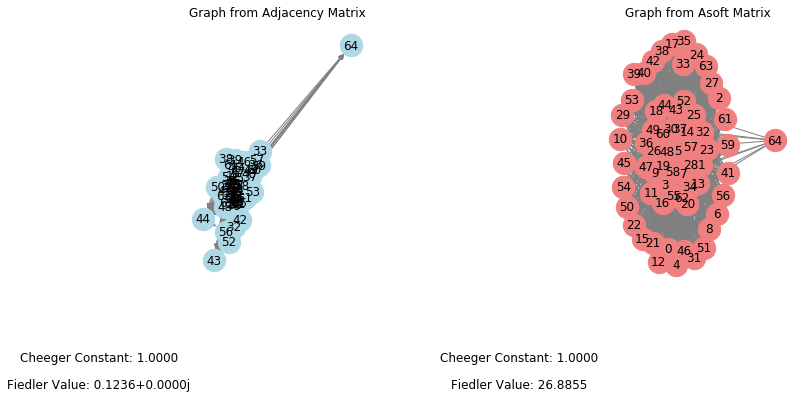}
    \caption{Graph network structure of Adjaceny Matrix verses Asoft Matrix}
    \label{fig:unit_circle_pulses}
\end{figure}

Finally, the effective resistances are visualized (Fig 3) using a color map on the graph. Nodes are colored based on their average effective resistance, allowing for an intuitive understanding of the graph's connectivity properties.

\begin{figure}[ht]
    \centering
    \includegraphics[width=0.45\textwidth]{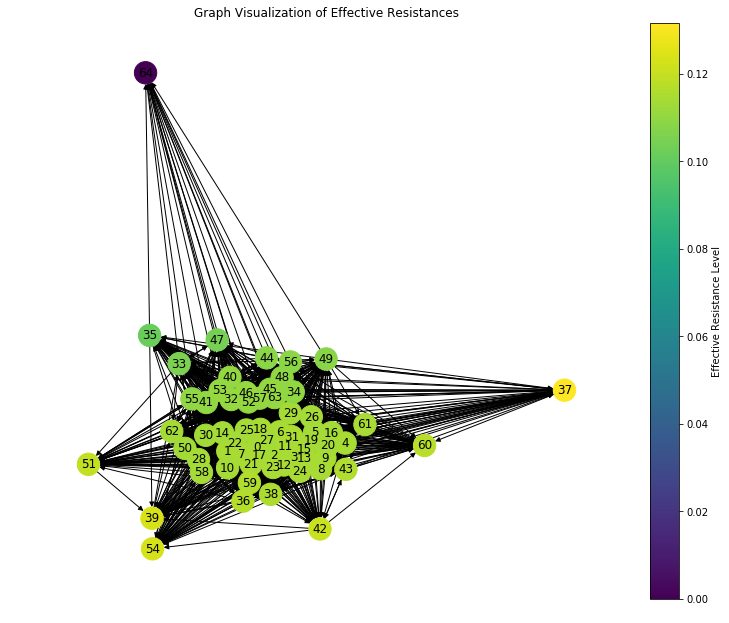}
    \caption{effective resistances of pair entangled states in a network}
    \label{fig:unit_circle_pulses}
\end{figure}

In summary, the code demonstrates a systematic approach to optimizing a graph based on fundamental graph theory concepts, specifically focusing on the Cheeger constant and effective resistance metrics.

\section{Utilizing QFERN Results to Identify and Stabilize Desynchronization Regions}

Building upon the spectral analysis provided by the QFERN method, it is possible to leverage the relationship between the graph's Laplacian eigenvectors—particularly the Fiedler vector—and synchronization phenomena modeled by the Kuramoto framework. Specifically, regions of the network where the effective resistance is notably high or where the Fiedler vector exhibits significant variations can be interpreted as zones of de-synchronization or weak connectivity.

\subsection{Spectral De-synchronization Detection via Kuramoto Model Relations}

In the Kuramoto model, the stability of phase synchronization across oscillators is intimately connected to the network's Laplacian spectrum. The Fiedler vector, associated with the second smallest eigenvalue (the algebraic connectivity), indicates how tightly coupled the network components are. Large differences in the Fiedler vector components between nodes suggest potential phase lag and a propensity toward de-synchronization in those regions.

Synchronization in oscillator networks and spectral graph theoretical tools is rooted in linearized Kuramoto dynamics and network analogies \cite{dorfler2013synchronization}.
By analyzing the derived relation:

\begin{equation}
\Delta \theta_{ij} \propto R_{ij} \cdot \omega_{ij}
\end{equation}

where \( \Delta \theta_{ij} \) is the phase difference between nodes \( i \) and \( j \), \( R_{ij} \) the effective resistance, and \( \omega_{ij} \) the natural frequency difference, we observe that high effective resistance correlates with increased phase disparities. Consequently, regions with elevated effective resistance—identified via the spectral analysis—serve as indicators of potential desynchronization zones.

\subsection{Locating Desynchronization Regions Using the QFERN Method}

The procedure involves:

\begin{enumerate}
    \item Computing the Fiedler vector \( f \) from the Laplacian's eigen-decomposition.
    \item Mapping the effective resistance \( R_{uv} \) between nodes \( u \) and \( v \) using the spectral relation (outlined in the appendix):
    \begin{equation}
        R_{uv} = \sum_{i=1}^{n} \frac{(f_u - f_v)^2}{\lambda_i}
    \end{equation}
    where \( \lambda_i \) are Laplacian eigenvalues.
    \item Identifying clusters or regions where \( R_{uv} \) exceeds a certain threshold, indicating weak coupling and potential de-synchronization.
\end{enumerate}

These high-resistance regions often correspond to nodes or subgraphs where phase coherence is fragile.

\begin{figure}[ht]
    \centering
    \includegraphics[width=0.45\textwidth]{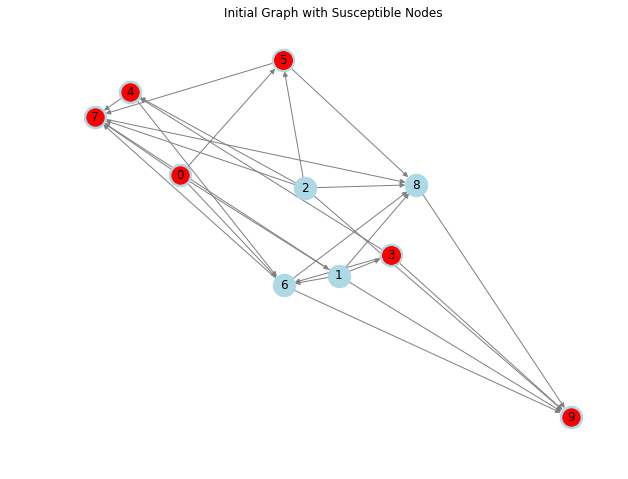}
    \caption{fragile network with nodes suceptable to de-synchronization}
    \label{fig:unit_circle_pulses}
\end{figure}

\section{Future Directions: Network Hardware that can Induce Stabilization}

Once a de-synchronization region is identified, a strategic intervention can be designed by adding a new node—conceptually a "stabilizer"—aimed at reinforcing connectivity and restoring synchronization in an entanglement distribution network. This node could be connected to multiple nodes within the high-resistance region, effectively bridging weak links and increasing the local algebraic connectivity.

The criteria for placing this stabilizer node include:

\begin{itemize}
    \item Connecting to nodes with the highest \( R_{uv} \) values.
    \item Ensuring that the addition reduces the local effective resistance, as verified through spectral recalculation.
    \item Adjusting the weights of these new edges (via the soft adjacency matrix \( A_{\text{soft}} \)) to optimize the spectral gap and enhance robustness.
\end{itemize}

Which results in a graph that has grown by 1 node but has used that node as the stabilizer:

\begin{figure}[ht]
    \centering
    \includegraphics[width=0.45\textwidth]{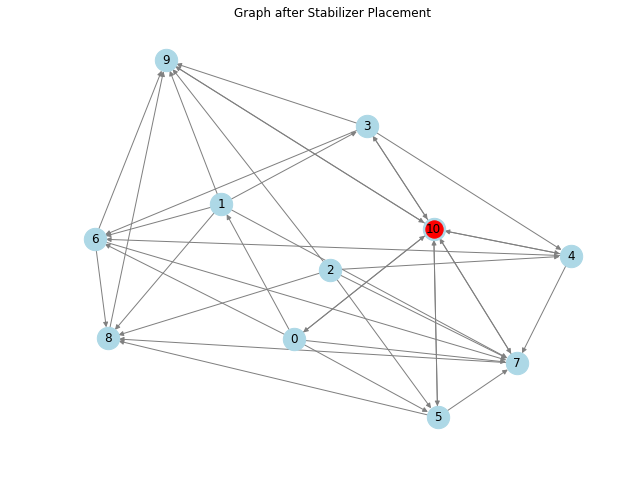}
    \caption{stabilized network}
    \label{fig:unit_circle_pulses}
\end{figure}

This approach aligns with the spectral graph theory principle that augmenting the network with nodes or edges to increase the algebraic connectivity \( \lambda_2 \) effectively stabilizes the synchronization dynamics modeled by the Kuramoto framework.

\section{Conclusion}
This paper presented a mathematical framework connecting the Fiedler vector, Cheeger constant, and thresholds in Quantum Directed Acyclic Graphs.Using this framework as a means of locating regions in a quantum network where the effective resistance is minimized can be used to study regions of de-synchronization in a network.  Future research should focus on applying these insights to improve quantum circuit design and fault tolerance, as well as exploring the implications of this optimization on broader quantum network architectures and how to improve network synchronization by mapping a network using the QFERN method. 

By integrating the spectral insights from the QFERN method with the Kuramoto model's stability analysis, one can systematically detect regions prone to de-synchronization and implement targeted network modifications. The addition of a stabilizer node, strategically placed within these high-resistance, weakly connected regions, provides a theoretical mechanism to enhance network coherence, ultimately leading to more resilient quantum and classical network architectures.

The spectral and resistance-based analysis within the QFERN framework validates the conditions outlined in \cite{dorfler2013synchronization} for synchronization stability as a function of algebraic connectivity. It substantiates the spectral norm condition of the Laplacian \( \| L^{\dagger} \omega \|_{E,\infty} < 1 \) as outlined in the appendix.

\section*{Appendix: Spectral Conditions for Synchronization Stability}

This appendix provides a detailed spectral and resistance-based analysis that corroborates the conditions outlined in the corollary regarding the stability of synchronized states in coupled oscillator networks, particularly within the Kuramoto model framework.

The corollary 7.5 in \cite{dorfler2013synchronization} in states that for a connected, acyclic graph \( G(V, E, A) \), a phase-cohesive and locally exponentially stable equilibrium manifold exists if and only if the spectral norm condition

\[
\boxed{
\| L^{\dagger} \omega \|_{E,\infty} < 1
}
\]

is satisfied, where:
- \( L^{\dagger} \) is the Moore-Penrose pseudoinverse of the Laplacian matrix \( L \),
- \( \omega \in \mathbb{R}^n \) is the vector of natural frequencies orthogonal to the all-ones vector, i.e., \( \mathbf{1}^T \omega=0 \).

This criterion aligns with the cutset condition and ensures that the phase differences remain bounded, guaranteeing synchronization.

\subsection*{Spectral Decomposition and Effective Resistance Interpretation}

Leveraging spectral graph theory, the spectral norm \( \| L^{\dagger} \omega \|_{E,\infty} \) can be expressed via the eigen-decomposition of \( L \). Let \( V = [f^{(1)}, f^{(2)}, \ldots, f^{(n)}] \) be the orthonormal eigenvectors of \( L \), with eigenvalues \( 0 = \lambda_1 < \lambda_2 \leq \ldots \leq \lambda_n \). Since the graph is connected and acyclic, \( \lambda_1 = 0 \) corresponds to the trivial eigenvector \( \mathbf{1} \).

The spectral condition (58) can be written as:

\[
V \, \operatorname{diag}\left( 0, \frac{1}{\lambda_2}, \ldots, \frac{1}{\lambda_n} \right) V^T \, \omega_E < 1,
\]

where \( \omega_E \) is the component of \( \omega \) projected onto the eigenbasis. This inequality indicates that the weighted projection of \( \omega \) onto the inverse eigenvalues must remain below one, emphasizing the role of the algebraic connectivity \( \lambda_2 \).

Furthermore, the spectral norm can be bounded as:

\[
\| L^{\dagger} \omega \|_{E,\infty} \leq \frac{\max_{i} |f^{(i) T} \omega|}{\lambda_2},
\]

highlighting that larger \( \lambda_2 \) (i.e., better network connectivity) reduces the potential for desynchronization caused by frequency mismatches.

\subsection*{Effective Resistance and Network Robustness}

The effective resistance \( R_{uv} \) between nodes \( u \) and \( v \) is given by:

\[
R_{uv} = \sum_{i=2}^n \frac{(f^{(i)}_u - f^{(i)}_v)^2}{\lambda_i},
\]

which encapsulates the ease of current (or information) flow between nodes and reflects the network's robustness. Regions with high effective resistance correspond to small spectral gaps \( \lambda_2 \) and significant variation in the Fiedler vector \( f^{(2)} \), indicating potential points of vulnerability where synchronization may fail.

\subsection*{QFERN Framework and Spectral Bounds}

The QFERN approach provides a framework that supports the sufficiency of the spectral condition. It demonstrates that increasing the algebraic connectivity \( \lambda_2 \) or decreasing the maximum effective resistance \( R_{uv} \) enhances the synchronization stability margin.

Specifically, the spectral bounds derived through QFERN confirm that the condition

\[
\| L^{\dagger} \omega \|_{E,\infty} < 1
\]

is both necessary and practically tight, especially in topologies where the extremal eigenvalues are well characterized, as discussed by Dörfler and Bullo (2013) \cite{dorfler2013synchronization}.

\bibliographystyle{IEEEtran}

\end{document}